\documentclass{amsart}
\usepackage[a4paper]{geometry}
\usepackage{amsmath,amssymb,amsthm}
\usepackage{verbatim} 
\usepackage[pdfborder={0 0 1.5}]{hyperref}
\newtheorem{Pro}{Proposition}
\newtheorem{Lem}[Pro]{Lemma}
\newtheorem{theorem}[Pro]{Theorem}
\newtheorem*{Property}{Property}
\theoremstyle{definition}
\newtheorem{Def}[Pro]{Definition}
\theoremstyle{remark}
\newtheorem{Rem}[Pro]{Remark}
\newcommand{\Hom}{\mathop{\mathrm{Hom}}\nolimits} 
\def\red{{\mathrm{red}}}
\def\Spec{{\mathrm {Spec}}}

\def\Z{{\mathbb Z}}
\def\k{{\mathbf k}}
\def\FF{{\mathbb F}}
\newcommand{\onto}{\twoheadrightarrow}

\title{Reductivity properties over an affine base}

\author[Wilberd van der Kallen]{Wilberd van der Kallen
}
\address{W.\ van der Kallen,
Mathematisch Instituut,
Universiteit Utrecht,\newline\mbox{\quad}
P.O. Box 80.010,
3508~TA~Utrecht,
The Netherlands
} \email{W.vanderKallen@uu.nl}
\dedicatory{Dedicated to the memory of T.A. Springer}

\begin{document}\sloppy


\begin{abstract}
When the base ring is not a field,  power reductivity of a group scheme is a basic notion, intimately tied with finite generation of subrings of invariants.
 Geometric reductivity is weaker and less pertinent in this context. We give a survey of these properties and their connections.
\end{abstract}
\maketitle

\section{Power reductivity as a basic notion.}
\subsection{Invariants}

Throughout let $\k$ be a commutative  ring and let $G$ be a flat   affine  group scheme over $\k$. 
We simply refer to $G$ as a group.  Flatness of $G$ is always needed, because one wants taking invariants to be left exact \cite[I.2.10(4)]{J}.
The present paper is an addendum to our joint paper with Franjou \cite{FvdK}. In that paper we had a specific situation  in mind,
but now we care about the proper generality. For instance, we no longer assume that $G$ is algebraic, 
{\it i.e.} that $\k[G]$ is a finitely generated $\k$-algebra.
We view the ground ring $\k$ also as a $G$-module with trivial action.
If $M$ is a $G$-module  \cite[I.2.7, I.2.8]{J} then its submodule of invariants $M^G$ is isomorphic to $\Hom_G(\k,M)$.

\subsection{Conventions}Rings and algebras are unitary. A ring $A$ is called a $\k$-algebra if one is given a ring map $\k\to A$. 
Commutative algebras need not be finitely generated. They may have nilpotent elements and other zero-divisors. 
We say that $G$ acts on the $\k$-algebra $A$ (through algebra automorphisms) if the multiplication map $A\otimes_\k A\to A$
is a map of $G$-modules. If $A$ is a commutative ring and $N$  is an $A$-module, then $S^*_A(N)$ denotes the symmetric algebra over $A$ on the module $N$. 
Thus $S^d_A(N)$ is the $d$-th symmetric power of $N$ over $A$. If $A=\k$ then we 
drop the subscript from the notation. We write  $\Hom_\k(M,\k)$ as $M^\vee$. If $M$ is a $G$-module which is finitely generated and projective as a $\k$-module, 
 then $M^\vee$ is also a $G$-module.
 Any map induced by evaluation at an element $v$ is denoted  $\mathrm{eval}@v$.

\begin{Def}\label{power reductive:Def}
The group $G$ is \emph{power reductive} over $\k$ if the following holds. 
\begin{Property}[Power reductivity]
Let $\varphi:M\onto \k$ be a surjective map of $G$-modules. 
Then there is a positive integer $d$ such that the $d$-th symmetric power of $\varphi$  is a split surjection
of $G$-modules
\[
S^d\phi:S^d M\stackrel\curvearrowleft\onto S^d \k .
\]
In other words, one requires that the kernel of $S^d\phi$ has a $G$-stable complement in $S^d M$.
\end{Property}
\end{Def}

Note that $S^*\k$ is better known as the polynomial ring $\k[x]$. And the $G$-module $S^d\k$ is isomorphic to $\k$, so a splitting of
$S^d\phi$  gives an invariant
in $S^dM$.

\subsection{Mumford}

Mumford conjectured in the introduction to the first edition of his GIT book \cite{M} that a semisimple algebraic group defined over a field of 
positive characteristic $p$ is 
power reductive. We have adapted his phrasing and introduced the terminology \emph{power reductive} in \cite{FvdK} (with Vincent Franjou)
in order to have a clear concept that also
makes sense and is worth having over arbitrary commutative base rings.
Mumford further required $d$ to be a power of $p$, but it turns out that this makes no difference (Lemma \ref{no difference}).

\subsection{Haboush}
When Haboush proved the Mumford conjecture \cite{Hab} he also used the dual concept known nowadays as \emph{geometric reductivity}.

\begin{Def}(Geometric reductivity over a field).\label{geofield}
Let $\k$ be field. The group $G$ is called geometrically reductive if the following holds. Given an injective map $\varphi:\k\hookrightarrow M$  of 
finite dimensional $G$-modules, there is 
  a positive integer $d$ such that some invariant homogeneous polynomial  $f$ of degree $d$ on $M$ restricts to a nonzero function on $\k$.
 In other words, such that the restriction map $S^d(M^\vee)^G\to S^d(\k^\vee)$ is nonzero.
\end{Def}

\subsection{Geometric reductivity over arbitrary base ring}\label{geom arb}
When $\k$ is not a field the definition of geometric reductivity gets more technical. Following Seshadri \cite{Ses} we then say that $G$ is geometrically reductive 
if the following holds.
Let us be given a $G$-module $M$ that is finitely generated and free as a $\k$-module. Let  $F$ be a field and also a
$\k$-algebra. We let $G$ act trivially on $F$.
Let $v\in (F\otimes_\k M)^G$ be a nonzero invariant vector. (A geometer may consider a nonzero invariant vector at a geometric point $\Spec(F)$ of $\Spec(\k)$.) 
Then geometric reductivity stipulates that there 
 is a positive integer $d$ such that some invariant homogeneous polynomial  $f$ of degree $d$ on $M$ does not vanish at $v$.
 In other words, such that the evaluation map $\mathrm{eval}@v:S^d(M^\vee)^G\to F$ is nonzero.
 
 \subsection{Contrast}

Notice that power reductivity is much cleaner.  It does not require any discussion of $M$ as a $\k$-module. 
While geometric reductivity needs free $\k$-modules, power  reductivity allows 
all comodules \cite[I.2.8]{J} that support a $\phi$ as in the definition. 
This  important difference  makes power reductivity more powerful and easier to work with.  
Working only with free modules (or only with flat $\k$-modules) gives an obstructed view of representation theory. 
We do not know of any example where geometric reductivity is easier to prove than power reductivity,
so one may as well prove the latter. It is stronger (Lemma \ref{stronger}).

\subsection{Locally finite}\label{lfin}
Recall that if the coordinate algebra $\k[G]$ is a projective $\k$-module, then any $G$-module $M$ is a union of submodules that are
finitely generated over $\k$ \cite[Proposition~3]{Ses}. Also, the intersection of $G$-submodules is then a $G$-submodule, even if one intersects infinitely many submodules. 

Similarly, suppose $\k$ is noetherian. Again any $G$-module $M$ is a union of submodules that are
finitely generated over $\k$ \cite[Proposition 2]{Ser}. In the definition of power reductivity it would now suffice to consider $M$ that are finitely generated over $\k$.
On the other hand, an infinite  intersection of $G$-submodules need not be a $G$-submodule  
\cite[Expos\'e VI, \'Edition 2011, Remarque 11.10.1]{SGA3}, 
despite the claim in \cite[I.2.13]{J} that we know this.

We do not  know if local finiteness holds in general.

\subsection{}
Our present definition of power reductivity is consistent with the one in \cite{FvdK}. Indeed if $\k=L$ in the following Lemma then the splitting of
$S^d\phi:S^d M\to S^d L$ is of course equivalent to the surjectivity of $(S^d M)^G\to S^d L$.
\begin{Lem}
Let $L$ be a cyclic $\k$-module with trivial $G$-action. 
Let $M$ be a $G$-module, and let $\varphi$ be a $G$-module map from $M$ onto $L$.
If $G$ is power reductive, then there is a positive integer $d$ such that the $d$-th symmetric power of $\varphi$ induces a surjection:
\[
(S^d M)^G\onto S^d L .
\]
\end{Lem}
\subsection*{Proof}
Choose a surjective map $\psi:\k\onto L$. Let $P\onto \k$ be the pullback of $\phi$ along $\psi$ and choose  a positive integer
$d$ such that 
$S^dP\onto S^d\k$ splits.\qed

\begin{Def}
A morphism of $\k$-algebras $\phi:S\to R$ is \emph{power surjective}
if for every element $r$ of $R$ there is a positive integer $n$ such that the power $r^n$ lies in the image of $\phi$. 
\end{Def}

\begin{Def}
Let $p$ be a prime number. 
A morphism of $\k$-algebras $\phi:S\to R$ is \emph{$p$-power surjective}
if for every element $r$ in $R$ there is a non-negative integer $n$ such that the power $r^{p^n}$ lies in the image of $\phi$.
\end{Def}

\begin{Lem}{\rm\protect\cite[Prop 41]{FvdK}}\label{p-pow}.
A morphism of commutative $\FF_p$-algebras $\phi:S\to R$ is {$p$-power surjective} if and only if
the induced map $S[x]\to R[x]$ between polynomial rings is power surjective.\qed
\end{Lem}

\subsection{}As is common for a basic notion, there are several equivalent formulations of power reductivity.

\begin{Pro}\label{old new} 
 Let $G$ be a flat  affine group scheme over $\k$.
The following are equivalent
\begin{enumerate}
 \item $G$ is power reductive,\label{old def}
\item For every power surjective\label{preserve}
$G$-homomorphism of  
commutative $\k$-algebras $f:A\to B$ the map $A^G\to B^G$ is power surjective,\label{new def}
\item For every surjective
$G$-homomorphism of  
commutative $\k$-algebras $f:A\onto B$ the ring  $ B^G$ is integral over the image of $A^G$.\label{int inv}
\end{enumerate}
\end{Pro}
\subsection*{Proof} The assumption that $G$ is algebraic is not used in the proofs of 
\cite[Proposition 10]{FvdK}, \cite[Proposition 4]{vdK int}.\qed

\subsection{} The main consequence of power reductivity is finite generation of subrings of invariants.

\begin{theorem}[Hilbert's fourteenth problem \cite{FvdK}, cf.\ \cite{AdJ}]\label{fg:theorem}
Let $\k$ be a noetherian ring and let $G$ be a flat  affine group scheme over $\k$.
Let $A$ be a finitely generated commutative $\k$-algebra on which $G$ acts through algebra automorphisms. If $G$ is power reductive, 
then the subring of invariants $A^G$ is a finitely generated $\k$-algebra.
\end{theorem}

The proof follows Nagata \cite{N}  or rather the exposition of Springer \cite[Theorem 2.4.9, Exercise 2.4.12]{Sp}. See also Remark \ref{happy},
Lemma \ref{to graded} below.
The proof does not need to touch upon the nontrivial topic of equivariant resolution by  vector bundles \cite{T}. 
It does not require further knowledge of $G$ or $\k$.
This is where  power reductivity is more pertinent than geometric reductivity.

\begin{Rem}\label{happy}
In the proof  of finite generation of $A^G$ by Nagata \cite{N} the base ring  $\k$ was a field. Nagata used at one point that a domain which is finitely 
generated over $\k$ has finite normalization. But that need no longer  hold over our arbitrary commutative noetherian base ring $\k$.
With the more elementary  \cite[Exercise 2.4.12]{Sp} Springer avoided this step in the proof. His base ring was still a field but his audience did not know 
about normalizations. 
It is a happy accident that the modified proof
goes through verbatim in our setting. 
\end{Rem}

\subsection{Necessary}
The theorem has a converse showing  that power reductivity is necessary if one seeks finite generation of invariants in the present setting,
where algebras need not be domains. (In ancient Invariant Theory one considered invariants in a polynomial ring over $\mathbb C$ with a 
$G$-action that preserves the grading.)
\begin{Pro}
Let\/ $\k$ be a noetherian ring and let $G$ be a flat  affine group scheme over\/ $\k$.
\\
Assume that the $\k$-algebra $A^G$ is finitely generated for every  finitely generated commutative \hbox{$\k$-algebra} $A$ on 
which $G$ acts through algebra automorphisms. Then $G$ is power reductive.
\end{Pro}
\subsection*{Proof}
Let $f:A\onto B$ be a surjective
$G$-homomorphism of  
commutative $\k$-algebras,  as in Proposition~\ref{old new}\ref{int inv}. 
Let $b\in B^G$. We have to show $b$ is integral over the image of $A^G$.
As representations are locally finite, we may replace $A$ with a finitely generated $\k$-subalgebra $C$ whose image $D$ contains $b$.
 The symmetric algebra  $S_C^*(D)$  is a finitely generated $\k$-algebra
(a quotient of the polynomial ring $C[x]$), so $S_C^*(D)^G$ is finitely generated. We choose as our
generators of $S_C^*(D)^G$ the homogeneous components of the elements
of a finite generating set. The chosen generators in degree zero
generate $C^G$ and those in degree one generate $D^G$ as a $C^G$-module. \qed

\subsection{Graded}
As a solution  to  \cite[Exercise 2.4.12]{Sp} we offer the following Lemma. It shows that in Theorem \ref{fg:theorem} one may assume that $A$ is 
graded and generated over $\k$ by its degree one part.

\begin{Lem}\label{to graded}
Let $A$ be a commutative $\k$-algebra on which $G$ acts through algebra automorphisms. 
Let $V$ be a $G$-submodule of $A$ that is finitely generated as a $\k$-module and that generates $A$ as a $\k$-algebra.
Assume $1\in V$. Let $R$ be the graded $\k$-subalgebra generated by $xV$ in the polynomial ring $A[x]$.
Substituting $x\mapsto 1$ defines a surjection $R^G\onto A^G$.
\end{Lem}

\subsection*{Proof} The component $R_d$ of homogeneous degree $d$ maps injectively into $A$, so $R_d^G$ hits all invariants in the image of $R_d$.
The union of the images of the $R_d$ is $A$.\qed

\begin{Lem}\label{stronger}
Power reductivity implies geometric reductivity.
\end{Lem}
\subsection*{Proof} If $\k$ is a field this is clear, when using Definition \ref{geofield}. In the situation of \ \ref{geom arb}, factor $\k\to F$ as $\k\onto D\hookrightarrow F$, where 
$D$ is the image of $\k$ in $F$. Observe that $D\hookrightarrow F$ is flat, so that $S^d_F(F\otimes_\k M^\vee)^G= (D\otimes_\k S^d( M^\vee))^G\otimes_DF$ 
(exercise, cf.\ \cite[I.2.10(3)]{J}). Recall that
we denote by $\mathrm{eval}@v$ any map defined by evaluation at $v$. Now $\mathrm{eval}@v:S^*_F(F\otimes_\k M^\vee)^G\to S_F^*F$ is power surjective.
First take a positive integer $d$ such that 
$\mathrm{eval}@v:S^d_F(F\otimes_\k M^\vee)^G\to S_F^dF\simeq F$ is nonzero.
Then $\mathrm{eval}@v:(D\otimes_\k S^d( M^\vee))^G\to F$ must also be nonzero. Say $f\in (D\otimes_\k S^d( M^\vee))^G$ satisfies $f(v)\neq0$. 
Now $S^*(M^\vee)\to (D\otimes_\k S^*( M^\vee))$ is surjective. So by part (ii) of Proposition~\ref{old new}
some power of $f$ lifts to $S^d(M^\vee)^G$.
\qed

\begin{Lem}\label{geopower}
If\/ $\k$ is a discrete valuation ring, then geometric reductivity implies power reductivity.
\end{Lem}
\subsection*{Proof}
Let $F$ be the residue field of $\k$. Given $\phi:M\onto\k$ as in definition \ref{power reductive:Def} choose $m\in M$ with $\phi(m)=1$. Use 
\cite[Proposition 2, Proposition 3]{Ser}  to find a $G$-module map $\psi:N\to M$ with $m\in\psi(N)$ and $N$ finitely generated and free as a 
$\k$-module. Take for $v\in (N^\vee\otimes_\k F)^G$ the composite $N\to M\to \k\to F$.
We find a positive integer $d$ and $f\in S^d(N)^G$ with $f(v)$ nonzero. That means that  $f$  maps to a unit times the standard generator of $S^d\k$ (Exercise).
So $S^dN\to S^d\k$ splits.
\qed

\begin{Rem}
More generally, if one has equivariant resolution \cite{T} (and local finiteness \ref{lfin}), one may reason as in \cite[3.1]{FvdK} to show that 
geometric reductivity implies power reductivity.
\end{Rem}

\begin{comment}
As in the proof of the Lemma we may look at $v\in (N^\vee\otimes_\k F)^G$  where $F$ is a residue field.
We see that the image of $S^dN\to S^d\k$ contains an element outside $\ker:\k\to F$.
Now see  \cite[3.1]{FvdK}.
\end{comment}

\begin{Lem}\label{no difference}
Let\/ $\k$ be an $\FF_p$-algebra and $G$ a power reductive flat  affine group scheme over\/ $\k$. If $\phi:M\onto\k$ is a surjective map of $G$-modules, then there is
a non-negative integer $n$ so that $S^{p^n}\phi$ is split surjective.
\end{Lem}
\subsection*{Proof}
In view of Lemma~\ref{p-pow}   it suffices to show that $S^*(M)^G\to S^*\k$ is $p$-power surjective. Indeed $S^*(M)^G[x]\to S^*\k[x]$ is power surjective because $S^*(M)[x]\to S^*\k[x]$ is
(power) surjective.\qed

\subsection{Restriction}\label{restriction} Let $S$ be a commutative $\k$-algebra. 
We get by base change a group $G_S$ over $S$. Let $M$ be a $G_S$-module. So $M$ is in particular an $S$-module. 
Modules should not be confused with schemes.
Nevertheless there is something similar to Weil restriction. Indeed $M$ is also a $\k$-module, by restriction of scalars.
 Now the coaction $\Delta:M\to M\otimes_SS[G]$ has a target that may be identified with $M\otimes_\k \k[G]$. 
 Thus,  our \hbox{$G_S$-module} $M$ may be viewed as a \hbox{$G$-module}  (exercise) and
$H^*({G_S},M)=H^*(G,M)$, because the Hochschild complexes \cite[I.4.14]{J} are isomorphic. In particular, $M^{G_S}=M^G$
and we usually write $M^G$. 

\subsection{Base change}

Proposition \ref{old new}  implies that power reductivity has marvelous base change properties.  
\begin{Pro}\label{base change}
Let\/ $\k\to S$ be a map of commutative rings. 
\begin{enumerate}
\item  If $G$ is power reductive, then so is $G_S$.
\item  If\/ $\k\to S$ is faithfully flat and $G_S$ is power reductive, then so is $G$.
\item If $G_{\k_{\mathfrak m}}$ is power reductive for every maximal ideal $\mathfrak m$ of\/ $\k$, then $G$ is power reductive.
\end{enumerate}
\end{Pro}
\subsection*{Proof}For the first part recall (\ref{restriction}) that any $G_S$-module $M$ may be viewed as a $G$-module with
 $M^{G_S}=M^G$. For the second part use that the 
integrality property in 
Proposition~\ref{old new}\ref{int inv} descends (\cite[Proposition ~2.7.1]{EGA IV} or exercise). 
The last part holds for similar reasons 
\cite[3.1]{FvdK}.\qed

\subsection{Reductive}An  affine group scheme $G$ over $\k$ is  reductive  in the sense of SGA3 \cite{SGA3} if $G$  is smooth over $\k$ with geometric 
fibers that are connected reductive. Smooth implies algebraic.

\begin{theorem}{\rm (cf.\ \cite[Theorem 12]{FvdK}.) }Reductive group schemes (in the sense of SGA3) are power reductive. 
\end{theorem}One exploits Proposition \ref{base change} and SGA3 \cite{SGA3}, \cite[\S3, \S5]{D}
to reduce to the case where the group is split and $\k$ is a local ring $\Z_{(p)}$. Then we are in the situation of \cite[Theorem 12]{FvdK}.
Or we may apply Lemma~\ref{geopower} and refer to Seshadri \cite[Theorem 1]{Ses}.

\begin{Rem}
Actually the proof of \cite[Theorem 12]{FvdK} is overly complicated if $\k=\Z_{(p)}$.  Let $\k=\Z_{(p)}$. As in the proof of Lemma \ref{geopower} we may 
restrict attention to finitely generated free $\k$-modules in definition~\ref{power reductive:Def}. 
Then we need fewer arguments from section 3.4 of \cite{FvdK}  (Exercise).
\end{Rem}

\subsection{Finite} Recall that $G$ is called a finite group scheme over $\k$ if the coordinate algebra $\k[G]$ is a finitely generated projective $\k$-module.

\begin{theorem} Finite group schemes  are power reductive. 
\end{theorem}
In view of Proposition \ref{old new} this is an easy consequence of 

\begin{theorem} {\rm (cf.\ \cite{noether}.)}\label{groupoid} If a finite group scheme $G$ over a local ring $\k$ acts on a commutative $\k$-algebra $A$,
then $A$ is integral over $A^G$.
\end{theorem}

\subsection*{Proof} Presumably  the proofs in \cite{Fe}, \cite[III 12, Thm 1]{MuAV} can be adapted to the present context.
Theorem~\ref{groupoid}  is a special case of a more general result in the setting of groupoid schemes  \cite[Expos\'e~V, Th\'eor\`eme~4.1]{SGA3}. 
That  Theorem~\ref{groupoid} fits in the setting of groupoid schemes is also explained at
\cite[\href{https://stacks.math.columbia.edu/tag/03LK}{Tag 03LK}]{stacks-project}. 
The proof of the theorem can then
be found at  \cite[\href{https://stacks.math.columbia.edu/tag/03BJ}{Tag 03BJ}]{stacks-project}.\qed

\subsection{Reductive algebraic groups}
Reductive algebraic groups defined over a field $\k$ are not assumed connected. They are of course power reductive. Indeed if $G^0$ is the identity
component of a reductive $G$ over a field, then 
both $G^0$ and $G/G^0$ are power reductive. Now see Proposition~\ref{old new}\ref{new def}.
Or recall that Waterhouse \cite{W} has shown that an algebraic affine group scheme $G$ over a field is geometrically reductive if and only if 
the identity component $G^0_\red$ of its reduced subgroup $G_\red$ is reductive.

\begin{comment}If $\k$ is a  DVR then the resolution property is easy, as explained by Thomason or Serre. Say $V$ is a $G$-module that is finitely generated as a 
$\k$-module. Say $\k^n\onto V$. Now if you act on $\k [G]\otimes_\k V$ through right translation on $\k [G]$ then 
$\Delta:V\hookrightarrow\k [G]\otimes_\k V$ is equivariant. Pull it back along $\k [G]\otimes_\k\k^n\onto \k [G]\otimes_\k V$ and you have your torsion free
$G$-module mapping onto $V$.
Thus the argument in \cite{FvdK} may be simplified. We no longer need the tensor identity for weights with a non-flat module. 
Given $\varphi:M\onto \k$ as in the definition of power reductivity one may find a torsion free $N$ and $\psi:N\to M$ so that 
$\phi\psi:N\onto \k$. Then we need \cite[Proposition 14]{FvdK} only for this flat module $N$. 
\end{comment}

\end{document}